\documentclass{wscpaperproc}

\usepackage{latexsym}
\usepackage{graphicx}
\usepackage{subcaption}
\usepackage{mathptmx}
\usepackage[T1]{fontenc}
\usepackage{booktabs}
\usepackage{multirow}
\graphicspath{{src/output/figures/}}

\usepackage{amsmath}
\usepackage{amsfonts}
\usepackage{amssymb}
\usepackage{amsbsy}
\usepackage{amsthm}

\usepackage[pdftex,colorlinks=true,urlcolor=blue,citecolor=black,anchorcolor=black,linkcolor=black]{hyperref}

\DeclareRobustCommand{\rev}[1]{\textcolor{black}{#1}}

\begin{document}

\pagestyle{fancyplain}

\thispagestyle{plain}
\firstPageHead{}

\chead{\fancyplain{}{\itshape Wang, Camur, \rev{Choudhuri,} and Li}}

\rhead{}
\cfoot{}
\renewcommand{\headrulewidth}{0pt} 

\makeatletter
\let\@internalcite\cite
\def\cite{\def\@citeseppen{-1000}%
    \def\@cite##1##2{(##1\if@tempswa , ##2\fi)}%
    \def\citeauthoryear##1##2##3{##1 ##3}\@internalcite}
\def\citeNP{\def\@citeseppen{-1000}%
    \def\@cite##1##2{##1\if@tempswa , ##2\fi}%
    \def\citeauthoryear##1##2##3{##1 ##3}\@internalcite}
\def\citeN{\def\@citeseppen{-1000}%
    \def\@cite##1##2{##1\if@tempswa, ##2)\else{}\fi}%
    \def\citeauthoryear##1##2##3{##1 (##3)}\@citedata}
\def\citeA{\def\@citeseppen{-1000}%
    \def\@cite##1##2{(##1\if@tempswa , ##2\fi)}%
    \def\citeauthoryear##1##2##3{##1}\@internalcite}
\def\citeANP{\def\@citeseppen{-1000}%
    \def\@cite##1##2{##1\if@tempswa , ##2\fi}%
    \def\citeauthoryear##1##2##3{##1}\@internalcite}
\def\shortcite{\def\@citeseppen{-1000}%
    \def\@cite##1##2{(##1\if@tempswa , ##2\fi)}%
    \def\citeauthoryear##1##2##3{##2 ##3}\@internalcite}
\def\shortciteNP{\def\@citeseppen{-1000}%
    \def\@cite##1##2{##1\if@tempswa , ##2\fi}%
    \def\citeauthoryear##1##2##3{##2 ##3}\@internalcite}
\def\shortciteN{\def\@citeseppen{-1000}%
    \def\@cite##1##2{##1\if@tempswa, ##2\else{}\fi}%
    \def\citeauthoryear##1##2##3{##2 (##3)}\@citedata}
\def\shortciteA{\def\@citeseppen{-1000}%
    \def\@cite##1##2{(##1\if@tempswa , ##2\fi)}%
    \def\citeauthoryear##1##2##3{##2}\@internalcite}
\def\shortciteANP{\def\@citeseppen{-1000}%
    \def\@cite##1##2{##1\if@tempswa , ##2\fi}%
    \def\citeauthoryear##1##2##3{##2}\@internalcite}
\def\citeyear{\def\@citeseppen{-1000}%
    \def\@cite##1##2{(##1\if@tempswa , ##2\fi)}%
    \def\citeauthoryear##1##2##3{##3}\@citedata}
\def\citeyearNP{\def\@citeseppen{-1000}%
    \def\@cite##1##2{##1\if@tempswa , ##2\fi}%
    \def\citeauthoryear##1##2##3{##3}\@citedata}
%
%
%
\def\@citedata{%
    \@ifnextchar [{\@tempswatrue\@citedatax}%
                  {\@tempswafalse\@citedatax[]}%
}

\def\@citedatax[#1]#2{%
\if@filesw\immediate\write\@auxout{\string\citation{#2}}\fi%
  \def\@citea{}\@cite{\@for\@citeb:=#2\do%
    {\@citea\def\@citea{, }\@ifundefined
       {b@\@citeb}{{\bf ?}%
       \@warning{Citation `\@citeb' on page \thepage \space undefined}}%
{\csname b@\@citeb\endcsname}}}{#1}}%

%
\def\@citex[#1]#2{%
\if@filesw\immediate\write\@auxout{\string\citation{#2}}\fi%
  \def\@citea{}\@cite{\@for\@citeb:=#2\do%
    {\@citea\def\@citea{; }\@ifundefined
       {b@\@citeb}{{\bf ?}%
       \@warning{Citation `\@citeb' on page \thepage \space undefined}}%
{\csname b@\@citeb\endcsname}}}{#1}}%

%
\def\@biblabel#1{}
\makeatother



\newdimen\bibindent
\bibindent=0.0em
\def\thebibliography#1{\section*{\refname}\list
   {}{\settowidth\labelwidth{[#1]}
   \leftmargin\parindent
   \itemindent -\parindent
   \listparindent \itemindent
   \itemsep 0pt
   \parsep 0pt}
   \def\newblock{}
   \sloppy
   \sfcode`\.=1000\relax}


\setlength{\baselineskip}{12.7pt}

\title{Data-Driven Stress Testing of Intermodal Freight Networks Using GAN-Generated Disruption Scenarios}

\author{\begin{center}Xudong Wang\textsuperscript{1}, Mustafa Can Camur\textsuperscript{1}, Sabarna Choudhuri\textsuperscript{2}, and Xueping Li\textsuperscript{1}\\
[11pt]
\textsuperscript{1}Dept.~of Industrial and Systems Eng., University of Tennessee, Knoxville, TN, USA\\
\textsuperscript{2}\rev{Amazon, New York, NY, USA} \end{center}
}

\maketitle

\vspace{-12pt}

\section*{ABSTRACT}
Intermodal freight networks are increasingly exposed to correlated, multi-mode disruptions, yet resilience assessments often rely on historical or uncorrelated scenarios that understate systemic risk. This paper develops a data-driven stress-testing framework integrating generative adversarial networks (GANs) with an intermodal optimization model to evaluate performance under realistic compound disruptions. \rev{The case study examines weather-related disruptions in the Tennessee Valley corridor.} Each GAN-generated scenario is used as a simulation input, and the resulting routing problem is solved to obtain system costs. Aggregating outcomes enables estimation of expected costs and identification of major risk drivers. Results show that historical disruptions increase total cost by about 3\%, whereas GAN-generated scenarios raise costs by over 25\%, producing an expected annual cost of \$5.11 million. \rev{Risk is concentrated in correlated multi-node failures and critical nodes such as the Port of Knoxville.} The framework helps identify vulnerabilities and prioritize resilience investments.

\section{INTRODUCTION}

\label{sec:introduction}

Freight disruptions pose a systemic risk to national-scale logistics networks that move massive volumes of goods across interconnected transportation modes. The 2022 Commodity Flow Survey \shortcite{USCensusCFS2025} reports that the U.S. businesses shipped 12.2 billion tons of goods valued at \$18.0 trillion. Trucks alone carried 68.1\% of total tonnage and 73.5\% of shipment value, indicating that highway capacity losses can rapidly propagate across intermodal supply chains. At the same time, infrastructure-disrupting extreme events are frequent. In 2024, the U.S. experienced 27 weather and climate disasters exceeding \$1B in damages \shortcite{NOAA_NCEI_Billions_US}, many of which disrupted intermodal transportation networks.

These conditions highlight the need for operationally meaningful stress testing of intermodal freight corridors. \rev{Recent freight-transport literature highlights growing attention to optimization, artificial intelligence, resilience, and disaster management \shortcite{KianiMavi2022FreightInnovation}; however, fewer studies explicitly focus on data-driven stress testing of correlated, multi-location disruption scenarios and their system-level impacts.} Transportation resilience research provides a complementary perspective, but it is typically divided between topology-based criticality measures and system-based consequence modeling \shortcite{MATTSSON201516}. Both approaches remain sensitive to the disruption scenarios, and tail events that drive worst-case losses are often underrepresented when scenario sets are small or manually constructed.


These limitations motivate scenario generation approaches that target realistic tail behavior. Rare-event simulation and adaptive stress testing provide principled mechanisms for exploring extreme trajectories, but may generate unrealistic stressors when domain constraints are insufficient \shortcite{BLANCHET201238,lee2015adaptive}. Deep generative models offer a complementary alternative: generative adversarial networks (GANs) can learn multivariate disruption structure from historical data, capturing spatial and temporal correlations without requiring explicit distributional assumptions \shortcite{Goodfellow2014GAN}. However, GAN-based approaches face challenges including training instability, limited rare-mode coverage, and the need to handle mixed discrete-continuous variables under operational constraints \shortcite{Gulrajani2017}.

This challenge is further amplified by the prevalence of correlated and compound disruptions. Research on interdependent systems shows that failures can propagate beyond their initial footprint and generate cascading effects across coupled networks \shortcite{Buldyrev2010}. From a hazards perspective, compound events arise from interacting drivers that jointly amplify impacts, with increasing evidence of correlated extreme weather generating concurrent infrastructure failures \shortcite{Zscheischler2018}, such as extreme weather at multiple locations. Yet many operational transport models treat disruptions as independent, which can significantly underestimate both the likelihood and severity of system-wide failures.

Risk-aware intermodal decision models address part of this gap through stochastic programming, robust optimization, and disruption-aware planning \shortcite{GBADEGOYE2025}. \rev{These models improve decision-making under uncertainty, but they still depend on predefined uncertainty sets, over-sampled scenarios, or historical cases.} When rare correlated disruptions are absent from those inputs, solutions can be brittle under compound conditions \shortcite{NEURIPS2019_8558cb40}. Poorly calibrated uncertainty sets risk excessive conservatism without capturing cross-modal correlations \shortcite{BenTalElGhaouiNemirovski2009}.

In this paper, we develop a data-driven stress-testing framework for intermodal freight networks. Historical disaster records are used to train a per-region GAN that generates correlated, multi-node disruption scenarios. Each scenario is mapped to disrupted network availability and evaluated with a corridor-scale intermodal optimization model, yielding a system cost defined as the sum of transportation cost, delay penalties, and emission cost across all shipments. Repeating this scenario-evaluation process over a large ensemble produces a distribution of outcomes that supports tail-risk quantification, expected cost estimation, and criticality ranking of infrastructure nodes. The main contributions are threefold: First, a data-driven GAN scenario generator that captures spatial and temporal disruption correlations; Second, a simulation framework that couples generative scenarios with optimization-based network evaluation; Third, a distributional stress-testing approach that quantifies expected loss and tail-risk concentrations under compound disruptions.

The remainder of the paper is organized as follows. Section~\ref{sect:2} formalizes the problem and presents the modeling framework for evaluating intermodal network performance under disruption scenarios. Section~\ref{sect:3} describes the methodology, including GAN-based scenario generation, the optimization model with the stress-testing procedure. Section~\ref{sect:4} presents the Tennessee Valley (TV) case study, including data sources, model calibration, and results under both historical and generated disruptions. Finally, Section~\ref{sect:5} concludes the paper and discusses implications for resilience planning and future research.

\section{PROBLEM STATEMENT}
\label{sect:2}
Intermodal freight networks operate as tightly coupled, multi-modal systems in which disruptions at key transfer hubs or river corridors can propagate across regions and transportation modes. In corridors such as the TV, freight flows depend heavily on river, rail, and highway interfaces concentrated at a limited number of nodes—a structural feature that creates a dual effect: localized failures may be absorbed through modal substitution when redundancy exists, but can trigger system-wide cost increases when rerouting options are constrained. Yet existing resilience assessment approaches largely rely on isolated historical events or deterministic worst-case assumptions, which fail to capture the correlated and compound nature of real-world disruptions. As a result, they cannot systematically evaluate how multi-node disruptions drive cost escalation or identify which infrastructure elements contribute most to systemic vulnerability.

\begin{figure}[htbp]
    \centering
    \includegraphics[width=\textwidth]{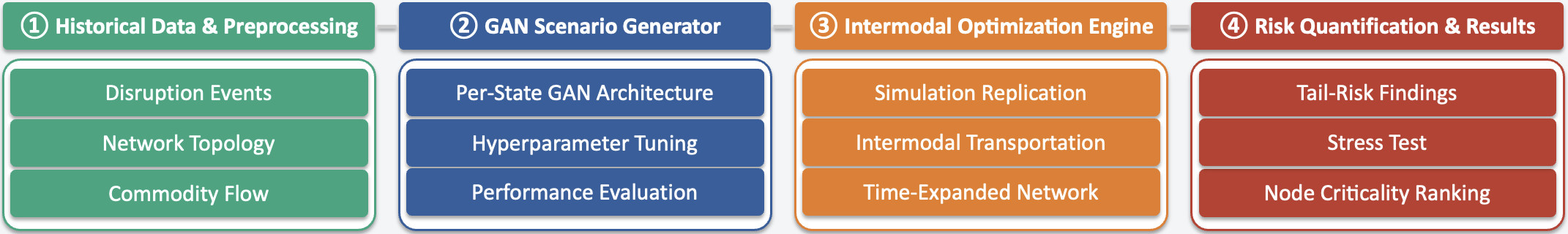}
    \caption{Data-driven stress-testing framework for intermodal freight networks, consisting of four stages: historical data assembly and preprocessing, GAN-based scenario generation, intermodal optimization-driven simulation, and risk quantification with node criticality ranking.}
    \label{fig:frameword}
\end{figure}

To address this limitation, we develop the data-driven stress-testing framework illustrated in Figure~\ref{fig:frameword}, enabling systematic evaluation of correlated disruption scenarios and their network-wide impacts. In this study, a disruption scenario is defined as a realization of simultaneous capacity reductions across a subset of network links or nodes over a horizon. The framework integrates four components: data assembly, GAN-based scenario generation, optimization-driven simulation, and risk quantification. Historical disruption records and network data are used to characterize system exposure. A GAN is then trained to generate realistic compound disruption scenarios that capture spatial and temporal dependencies beyond the historical record. Each scenario is evaluated through a time-expanded intermodal optimization model by applying the corresponding capacity reductions and re-optimizing flows, and its impact is quantified as the increase in system cost relative to a baseline under nominal network conditions. Aggregating these cost increases across scenarios produces a distribution of outcomes that enables expected cost estimation, tail-risk identification, and node criticality ranking, supporting data-driven identification of vulnerable and critical infrastructure.

\section{METHODOLOGY}
\label{sect:3}
\begin{figure}[htbp]
        \centering
        \includegraphics[width=\textwidth]{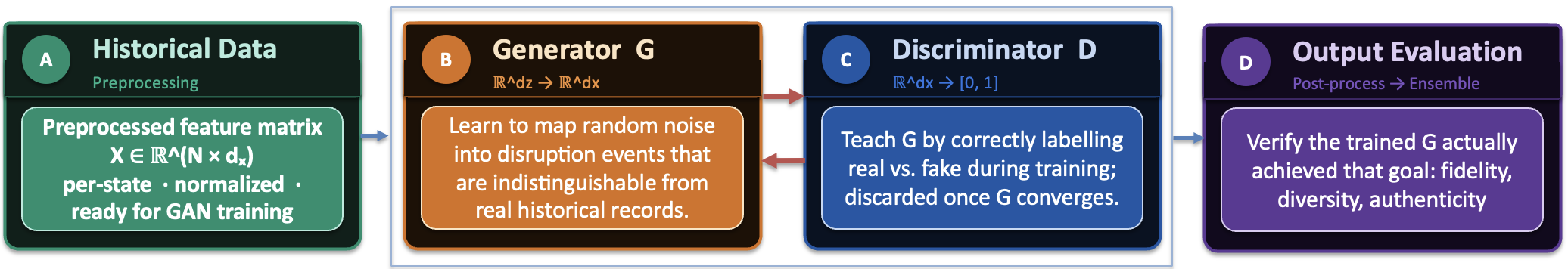}
        \caption{The framework of the GAN in our experiment.}
        \label{fig:GAN_framework}

\end{figure}

The proposed framework consists of two methodological components. The first is a data-driven scenario engine: a stratified, per-region GAN trained on historical event data that generates synthetic disruption scenarios. The regional decomposition is implemented at the state level, selected based on empirical evaluation of alternative spatial granularities. A single GAN trained on the full network produced diffuse scenarios that failed to capture localized disruptions, while a finer county-level decomposition led to overfitting and limited diversity. The state-level formulation provides a balance between spatial fidelity and statistical robustness. The second component evaluates each generated scenario by solving the intermodal routing model under the corresponding disrupted network conditions, producing one system cost observation per scenario. \rev{This creates a stress-testing experiment in which the scenario inputs are learned from historical disruption data rather than sampled from independent assumptions.} Figure \ref{fig:GAN_framework} shows the framework of the GAN in our experiment, and we will explain the detailed steps in the following subsections.

\subsection{GAN-Based Disruption Scenario Generation}
\label{sect:3.1}
The first component of the proposed framework is a data-driven scenario engine designed to synthesize realistic, multi-node disruption events, where each scenario specifies simultaneous capacity reductions across a subset of network links or nodes over time. By utilizing a stratified, per-region GAN architecture trained on historical disaster records, the engine captures spatial and temporal correlations that traditional independent models often overlook, enabling the generation of compound and tail-risk disruption patterns that drive systemic cost increases within the intermodal network. \rev{Here we use administrative regions because disruption records are reported with county-level identifiers and because county-level outputs can be directly linked to resource allocation, while the same framework can accommodate regions based on corridors, hazards, or infrastructures when such partitions are available.}

\subsubsection{Feature Representation and Preprocessing}
Let $\mathcal{R}$ denote the set of geographic subregions used to stratify the historical data. Each disruption event occurring in subregion $r \in \mathcal{R}$ is represented by a multivariate feature vector $\mathbf{x}^{(r)} \in \mathbb{R}^{d_x}$, where $d_x$ is the number of feature dimensions capturing temporal, spatial, categorical, and impact-related attributes. Temporal variables such as time-of-year and time-of-day are encoded using sine–cosine pairs to preserve their cyclical structure, ensuring that adjacent periods across year or day boundaries remain proximate in feature space. Spatial location is represented by the subregion index $r$, and event type is encoded as an integer category. The feature vector $\mathbf{x}^{(r)}$ provides a structured representation of each disruption event, and the variables described above illustrate the feature construction rather than define an exhaustive set.

Impact-related continuous variables in $\mathbf{x}^{(r)}$ often exhibit heavy-tailed distributions that can hinder stable training of generative models. For each such component $x_j^{(r)}$, we apply the log transform $x_j^{(r)\prime} = \log(1 + x_j^{(r)})$ to reduce skewness, followed by feature-wise standardization to zero mean and unit variance so that all variables operate on comparable scales during training \shortcite{Xu2019CTGAN}. The resulting normalized vector, denoted by $\tilde{\mathbf{x}}^{(r)}$, defines a homogeneous feature space for GAN learning, while the invertibility of these transformations allows generated samples to be mapped back to their original physical units after synthesis.

\subsubsection{Per-Region GAN Architecture and Training} \label{Per-Region GAN Architecture and Training}
Rather than fitting a single global generative model, we train one independent GAN per subregion $r \in \mathcal{R}$. This stratified design allows each generator to specialize in the local hazard profile and spatial distribution of its region, avoiding the distortion that would arise if regions with higher disruption event frequencies dominate the gradient signal in a pooled training objective.

Each regional GAN consists of two adversarially trained neural networks. The generator $G_{\theta_G}^{(r)}: \mathbb{R}^{d_z} \to \mathbb{R}^{d_x}$, parameterized by $\theta_G$, maps a latent noise vector $\mathbf{z} \sim \mathcal{N}(0, I_{d_z})$ of dimension $d_z$ to a synthetic disruption vector $\hat{\mathbf{x}} = G_{\theta_G}^{(r)}(\mathbf{z}) \in \mathbb{R}^{d_x}$. The generator uses fully connected layers with batch normalization and ReLU activations and a linear output layer. The discriminator $D_{\theta_D}^{(r)}: \mathbb{R}^{d_x} \to [0,1]$, parameterized by $\theta_D$, mirrors this depth using LeakyReLU activations and a sigmoid output.

Training minimizes the standard adversarial objective \shortcite{Goodfellow2014GAN}:
\begin{equation}
\min_{\theta_G} \max_{\theta_D} \; \mathbb{E}_{\mathbf{x} \sim p_{\text{data}}^{(r)}} \big[ \log D_{\theta_D}^{(r)}(\mathbf{x}) \big] + \mathbb{E}_{\mathbf{z} \sim \mathcal{N}(0,I_{d_z})} \big[ \log (1 - D_{\theta_D}^{(r)}(G_{\theta_G}^{(r)}(\mathbf{z}))) \big],
\end{equation}
where $p_{\text{data}}^{(r)}$ is the empirical distribution of preprocessed historical disruptions for subregion $r$. Both networks are optimized alternately within each mini-batch using the Adam optimizer.

Network architecture and training hyperparameters are selected via a grid search. Each candidate configuration is scored using a three-axis quality framework adapted from precision, recall, density, and coverage metrics \shortcite{kynkaanniemi2019improved}: \emph{fidelity} quantifies distributional similarity between real and synthetic samples using mean column-wise Wasserstein-1 distance as an evaluation metric \shortcite{Gulrajani2017}; \emph{diversity} measures the fraction of real events covered by the synthetic support; and \emph{authenticity} measures the fraction of synthetic events that fall within the real data manifold. The composite score is the unweighted mean of the three axes, and the configuration achieving the highest score is retained.

\subsubsection{Synthetic Scenario Generation and Post-Processing}
Once trained, each regional generator $G_{\theta_G}^{(r)}$ serves as a stochastic scenario sampler. Independent noise vectors $\mathbf{z} \sim \mathcal{N}(0, I_{d_z})$ are passed through $G_{\theta_G}^{(r)}$ to produce normalized synthetic vectors $\hat{\mathbf{x}} \in \mathbb{R}^{d_x}$. Each component $\hat{x}_j'$ is then inverted: for heavy-tailed components the inversion is $\hat{x}_j = e^{\hat{x}_j'} - 1$; for all components, inverse standardization restores original scale. Negative values are clipped to zero.

Discrete attributes require additional post-processing: subregion and category identifiers are rounded and clipped to valid ranges, then mapped back to their labels using encoding dictionaries constructed during preprocessing. The per-region synthetic datasets are combined into an ensemble of disruption scenarios that serve as simulation inputs to the intermodal optimization model across the stress-testing experiment.

\subsection{Intermodal Transportation Optimization and Stress Testing}

Each GAN-generated disruption scenario specifies a set of affected geographic locations and a disruption time window. For each scenario, the affected locations are translated into network availability: any edge whose route traverses an affected location during that window is removed from the feasible path set. The intermodal optimization model is then solved under this modified network, yielding a total system cost observation. Executing this process across the full scenario ensemble produces a distribution of cost outcomes that supports tail-risk quantification, expected loss estimation, and criticality analysis of individual infrastructure nodes. \rev{In this way, the optimization model evaluates the network response for each scenario, while the simulation experiment is formed by repeated scenario evaluation and aggregation.}

\subsubsection{Model Formulation}

This model serves as the simulation engine, solved once per replication; when a disruption scenario is applied, edge availability is modified and the model re-solved to yield that replication's system cost.

Let $N$ denote the set of nodes, $E$ the set of directed edges, $C$ the set of containers, $T$ the set of discrete time periods, and $\mathcal{M}^E$ the set of transportation modes including waiting at a node. The subset $\mathcal{M}^F \subset \mathcal{M}^E$ contains full-container-load (FCL) modes, which consolidate multiple containers into scheduled vehicles. For each container $c \in C$, let $\mathcal{P}_c$ be its set of feasible paths from source node $v_{cs}$ to demand node $v_{cd}$, and for each path $p \in \mathcal{P}_c$, let $\mathcal{O}_p$ be the set of transportation method combinations. The earliest available time of container $c$ is $a_c$ and its delivery deadline is $\nu_c$.

Let $\tau_{povmt}$, the travel time to reach node $v$ via mode $m$ under path $p$, combination $o$, starting at time $t$; $\epsilon_{po}$, the transportation cost of path $p$ under combination $o$; $\kappa_e$, the carbon emission cap on edge $e$; $\phi_{em}$, the carbon emissions per container traversing edge $e$ by mode $m$; $\lambda_e$, the edge container capacity per period; $\sigma_m$, the FCL vehicle container capacity; $\theta_{em}$, the maximum number of available vehicles on edge $e$ by mode $m$ per period; $\alpha$, the carbon tax per unit of excess emission; and $\beta$, the unit lateness penalty.

The decision variables are: $X_{cpt} \in \{0,1\}$, whether container $c$ takes path $p$ departing at time $t$; $Y_{cpo} \in \{0,1\}$, whether it uses transportation combination $o$ on path $p$; $V_c \in \mathbb{Z}_+$, the arrival time of container $c$; $L_c \in \mathbb{Z}_+$, its lateness beyond deadline; and $Q_e \in \mathbb{R}_{\geq 0}$, excess carbon emission on edge $e$. Auxiliary binary variables $A_{cimt}$, $B_{cimt}$, and $F_{cemt}$ track container location state and edge occupancy at each time period.

The objective minimizes total transportation cost, lateness penalties, and carbon emission tax:
\begin{equation}
\min \;\; \sum_{c \in C}\!\left(\beta L_c + \sum_{p \in \mathcal{P}_c}\sum_{o \in \mathcal{O}_p} \epsilon_{po}\, Y_{cpo}\right) + \alpha \sum_{e \in E} Q_e. \label{eq:obj}
\end{equation}

The key constraints are as follows. Each container selects exactly one path departing within its availability window, and exactly one transportation combination for that path:
\begin{equation}
\sum_{p \in \mathcal{P}_c}\sum_{t \geq a_c} X_{cpt} = 1, \quad \sum_{o \in \mathcal{O}_p} Y_{cpo} = \sum_{t \geq a_c} X_{cpt} \quad \forall\, c, p.
\end{equation}

Arrival time and lateness are computed from the selected path, combination, and departure time:
\begin{equation}
\sum_{p, o} \tau_{pov_{cd}m_o t}\,Y_{cpo} + t\,X_{cpt} \leq V_c, \quad L_c \geq V_c - \nu_c \quad \forall\, c, t.
\end{equation}

Binary state variables $A_{cimt}$ represent that container $c$ arrives at node $i$ via mode $m$ at time $t$, and $B_{cimt}$ that container $c$ is in transit toward or waiting at node $i$ via mode $m$ at time $t$. Variable $F_{cemt}$ indicates edge occupancy. These variables satisfy mutual exclusivity, ensure consistency between arrivals and departures, and enforce flow continuity along edges; details are omitted for brevity.

Carbon emissions on each edge are limited, with overages taxed:
\begin{equation}
\sum_{c \in C} \phi_{em}\, F_{cemt} \leq \kappa_e + Q_e \quad \forall\, e, m, t.
\end{equation}

Edge flow capacity, vehicle availability, and FCL consolidation capacity are enforced:
\begin{equation}
\sum_{c}\sum_{m} F_{cemt} \leq \lambda_e, \quad
\sum_{c} F_{cemt} \leq \theta_{em}, \quad
\sum_{c} F_{cemt} \leq \sigma_m \;\;(m \in \mathcal{M}^F) \quad \forall\, e, t.
\end{equation}

\subsubsection{Disruption Integration and Stress-Testing}

A disruption scenario $\boldsymbol{\xi}$ from the GAN specifies a set of affected locations and a time window $[t_s, t_e]$. Any edge $e \in E$ that passes through an affected location is made unavailable during $[t_s, t_e]$ by removing it from the feasible path set: for all paths $p$ that traverse $e$ and all departure times $t$ whose travel window overlaps $[t_s, t_e]$, the corresponding $X_{cpt}$ variables are fixed to zero. Compound scenarios disable multiple edges simultaneously, reflecting concurrent multi-location events drawn from the GAN ensemble.

For each scenario evaluation, the optimization model is re-solved under the disrupted network. Let $C(\boldsymbol{\xi})$ denote the optimal cost under scenario $\boldsymbol{\xi}$ and $C_0$ the baseline under nominal availability; the scenario cost increase is $\Delta C(\boldsymbol{\xi}) = C(\boldsymbol{\xi}) - C_0$. Aggregating across $n$ generated scenarios yields an empirical distribution of resilience impacts; the probability-weighted expected loss is \rev{\(\mathbb{E}[\Delta C] = \sum_{i=1}^{n} \pi_i \,\Delta C(\boldsymbol{\xi}_i)\)}, where $\pi_i$ is the probability weight of scenario $\boldsymbol{\xi}_i$ estimated from GAN output frequencies. Node criticality is measured as the expected $\Delta C$ over all scenario evaluations in which that location is disrupted, providing a scenario-based ranking of infrastructure importance.

\section{CASE STUDY: TV REGION INTERMODAL NETWORK}
\label{sect:4}
We apply the proposed framework to the TV transportation corridor, a critical freight network serving the southeastern U.S. This case study demonstrates how adversarial scenario generation identifies vulnerabilities in real-world infrastructure systems and informs resilience investments for the Tennessee Valley Authority (TVA) and regional transportation agencies.

\subsection{Data Sources}
\begin{figure}[htbp]
    \centering
    
    \begin{minipage}[b]{0.57\textwidth}
        \centering
        \includegraphics[width=\textwidth]{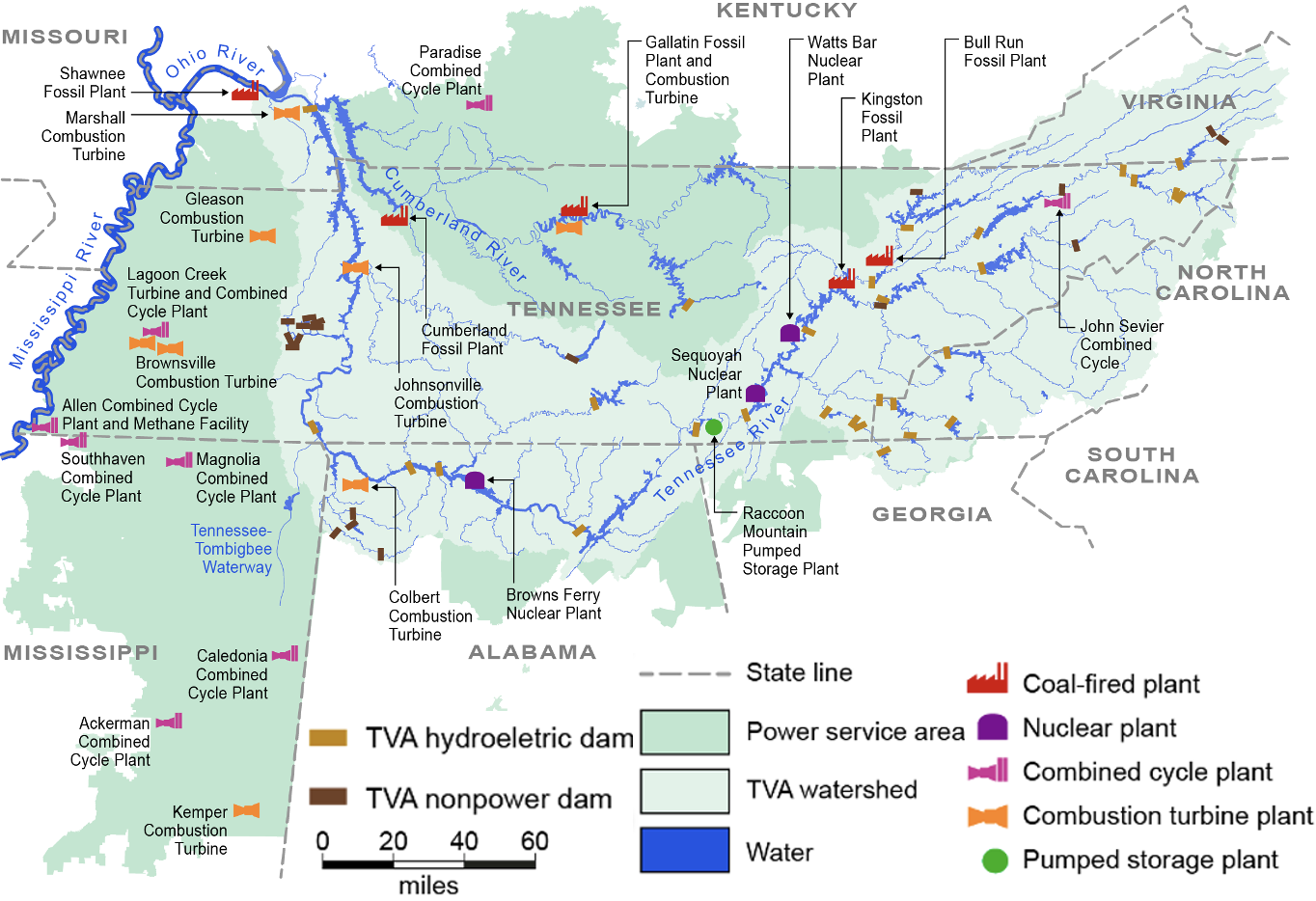}
        \caption{TVA's service area and facilities.}
        \label{fig:region}
    \end{minipage}
    \hfill
    \begin{minipage}[b]{0.42\textwidth}
        \centering
        \includegraphics[width=\textwidth]{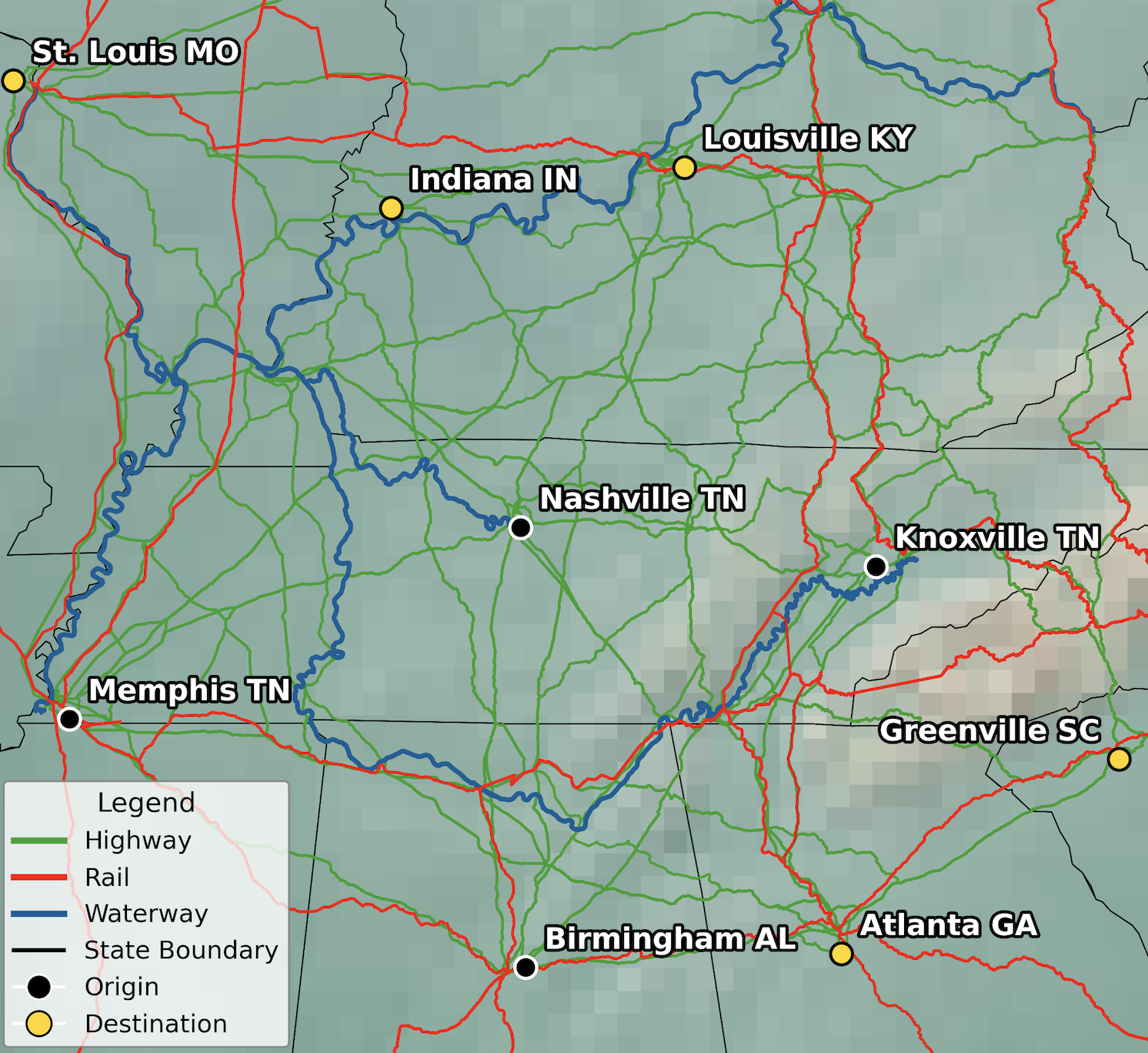}
        \caption{Multimodal freight network with OD pairs in the TVA region.}
        \label{fig:TVA}
    \end{minipage}

\end{figure}

The TV intermodal corridor connects major metropolitan areas including Memphis, Nashville, Chattanooga, and Huntsville, and links to national gateways like Savannah and Norfolk. Figure~\ref{fig:region} shows the geographic extent of the TVA service region along with key infrastructure assets \shortcite{gao2022tva}.

Highway, rail, and inland waterway network geometries are obtained from the Freight and Fuel Transportation Optimization Tool \cite{FTOT}, which provides detailed multimodal route shape files for roadway, railway, and waterway networks. In this study, the FTOT network serves as the unified infrastructure layer for constructing the intermodal network. Intercity freight demand is calibrated using the Freight Analysis Framework, which reports commodity-level origin–destination (OD) flows between metropolitan areas \cite{FAF5}. These datasets are integrated to form the multimodal transportation network in the TVA region, as shown in Figure~\ref{fig:TVA}.

Disruption data are drawn from the NCEI Storm Events Database \cite{NCEI_StormData}, which records county-level severe weather events across the U.S. For TVA counties, we extract event timing, geographic identifiers, hazard types, and impact metrics. The resulting dataset characterizes spatial and temporal patterns of disruptions and provides the basis for scenario generation and stress-testing analysis.

\subsection{Disruption Case Study}

This case study evaluates the proposed framework on a real-world intermodal freight corridor in the TV region, which is characterized by strong dependence on inland waterways, rail connections, and highway interfaces. The analysis aims to assess how correlated disruption scenarios affect network performance and to identify infrastructure components that contribute most to systemic vulnerability. Using historical disaster data and a detailed corridor-scale network model, we generate and evaluate a large ensemble of disruption scenarios to quantify cost impacts and characterize resilience under realistic stress conditions.

\subsubsection{GAN-Generated Disruption Scenarios}

We first analyze historical severe weather events extracted from the NCEI Storm Events Database. County-level annual disruption frequencies, hazard types, and severity measures are aggregated to characterize spatial exposure across the TV. To extend disruption modeling beyond the finite historical record, we train a per-state GAN on structured county-level disaster data for the five TVA states. Each generator learns the joint distribution of spatial location, hazard characteristics, and annual event intensity, enabling the synthesis of realistic disruption configurations.

The trained GAN produces synthetic annual disruption scenarios that preserve multivariate dependencies across counties and hazard types. Compared to historical observations, generated scenarios exhibit similar spatial clustering patterns while allowing moderate upper-tail variation and increased multi-county concurrence, thereby capturing plausible compound disruption structures.

\begin{figure}[!ht]
    \centering
    
    \begin{subfigure}[b]{0.48\textwidth}
        \centering
        \includegraphics[width=\textwidth]{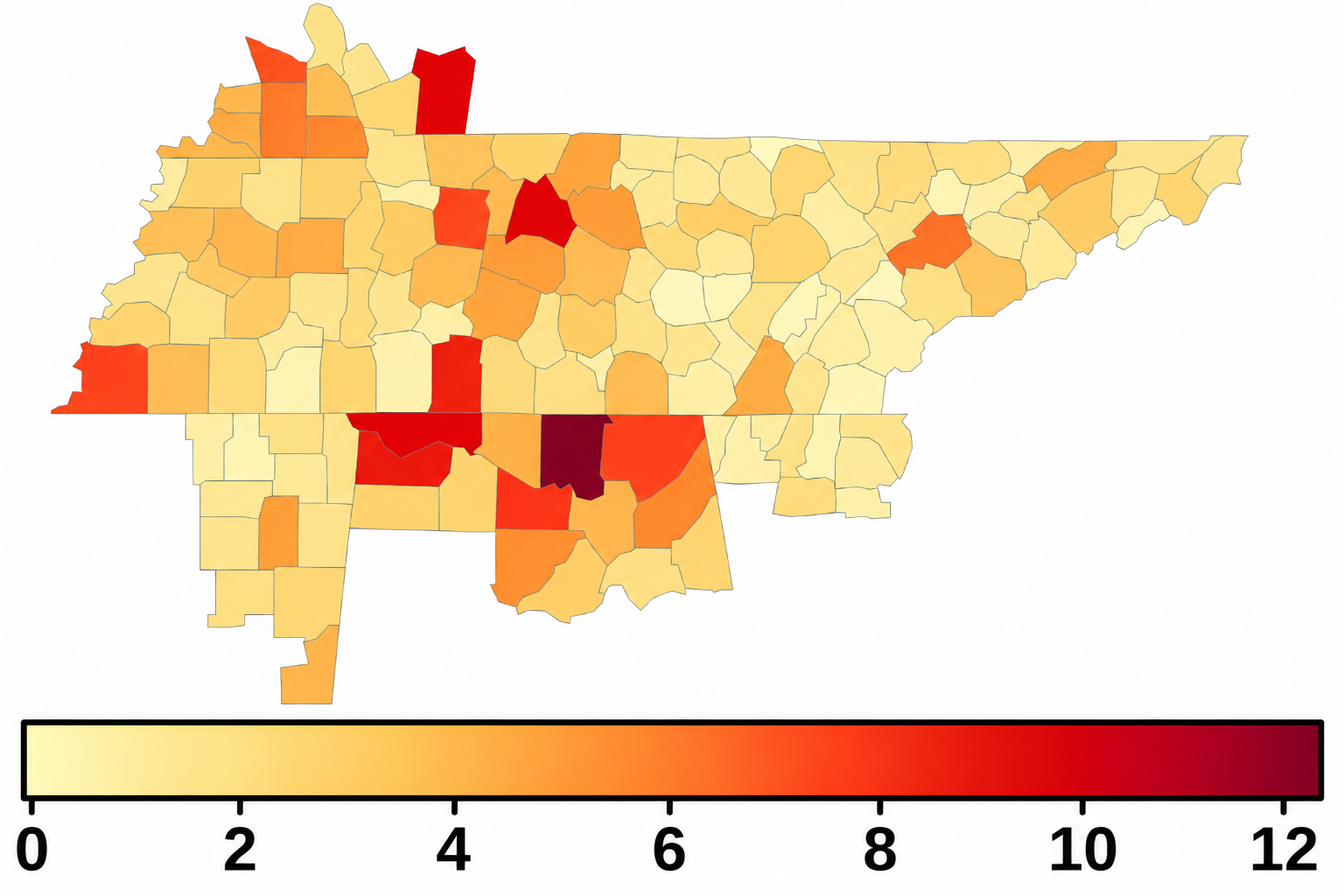}
        \caption{Average annual disruption cases since 2015.}
        \label{fig:data}
    \end{subfigure}
    \hfill
    \begin{subfigure}[b]{0.48\textwidth}
        \centering
        \includegraphics[width=\textwidth]{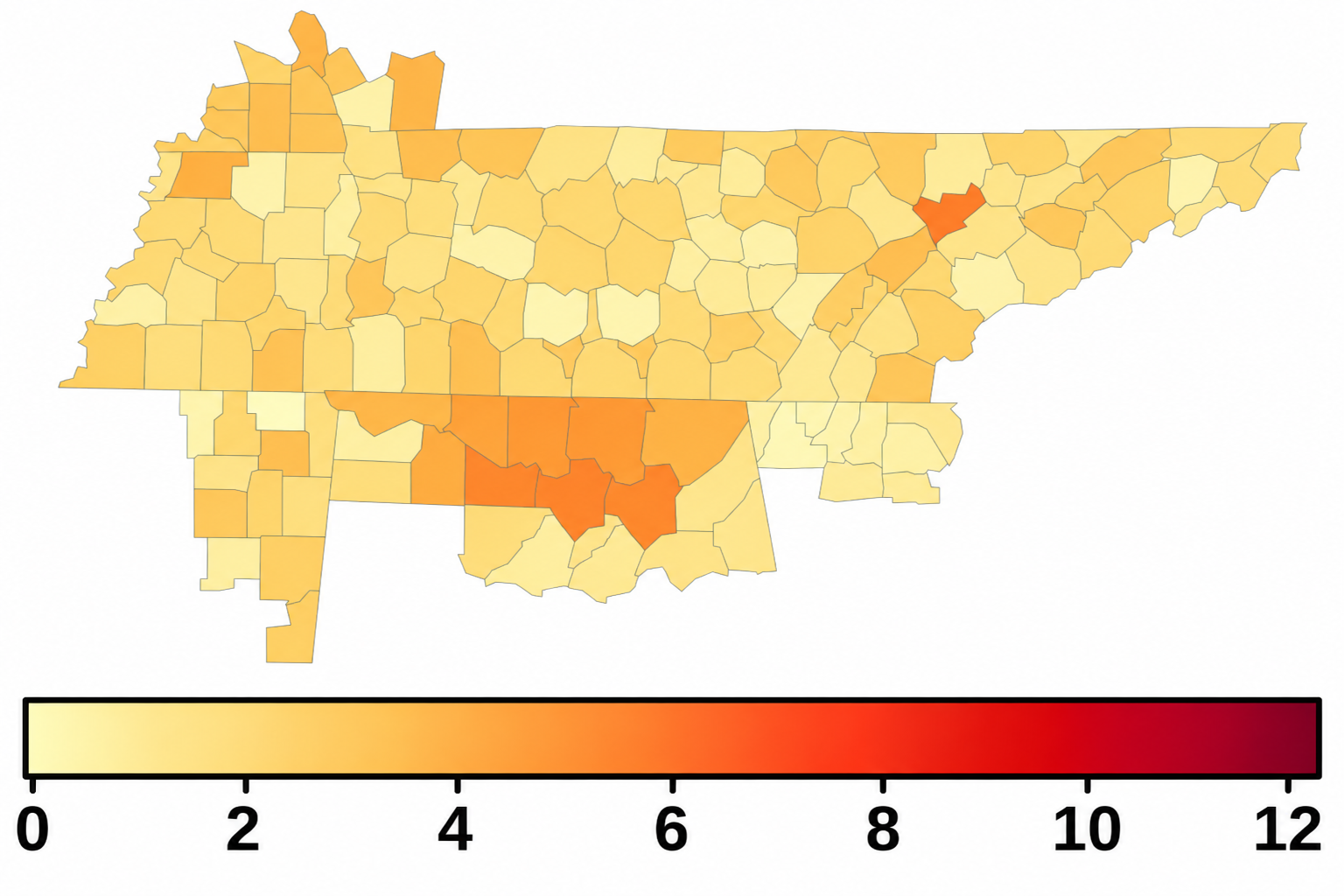}
        \caption{GAN-generated average annual disruption cases.}
        \label{fig:gan}
    \end{subfigure}
    
    \caption{Comparison of county-level annual disaster frequency: historical record in the last 10 years versus GAN-generated synthetic scenarios.}
    \label{fig:tva_freq_comparison}
\end{figure}

Figure~\ref{fig:tva_freq_comparison} compares the historical 10-year disruption frequency distribution with synthetic output from the best-performing configuration. The GAN preserves the rank ordering of high-risk regions, such as the Memphis metropolitan area, the Nashville basin, and northeastern Tennessee, while smoothing isolated extremes. \rev{This smoothing helps the generator learn broader spatial patterns from sparse historical records, and the resulting scenarios are used as plausible stress inputs for downstream tail-risk and node-criticality analysis.}

\subsubsection{GAN Evaluation and Performance}

\rev{Using the GAN performance indices introduced in Section \ref{Per-Region GAN Architecture and Training}, we select the best-performing configuration as the one that maximizes the composite score over fidelity, diversity, and authenticity. A grid search over 81 hyperparameter combinations is conducted, varying learning rate $\eta \in \{10^{-3},\,5\times10^{-4},\,10^{-4}\}$, latent dimension $d_z \in \{16,\,32,\,64\}$, hidden width $d_h \in \{128,\,256,\,512\}$, and training epochs $E \in \{300,\,500,\,800\}$. One GAN is trained per state, and scores are macro-averaged across states to determine the final configuration.}

Across all trials, fidelity remains high, ranging from 0.827 to 0.906 with a mean of 0.871, indicating stable marginal distribution learning. Diversity varies more, ranging from 0.589 to 0.723 with a mean of 0.682, and serves as the primary factor distinguishing stronger from weaker models due to the difficulty of capturing infrequent event modes. Authenticity shows moderate variation with a mean of 0.786, reflecting a balance between realism and avoidance of memorization. The composite score ranges from 0.701 to 0.805 with a mean of 0.780. The best-performing configuration achieves 0.805 with $\eta = 10^{-4}$, $d_z = 16$, $d_h = 512$, and $E = 800$, and its per-state performance is summarized in Table~\ref{tab:best_state}.




\begin{table}[!ht]
\centering
\caption{Per-state metrics for best GAN configuration.}
\label{tab:best_state}
\begin{tabular}{lcccc}
\toprule
State & Fidelity & Diversity & Authenticity & Composite \\
\midrule
Alabama (AL)     & 0.910 & 0.706 & 0.776 & 0.797 \\
Georgia (GA)     & 0.849 & 0.643 & 0.800 & 0.764 \\
Kentucky (KY)    & 0.889 & 0.740 & 0.836 & 0.822 \\
Mississippi (MS) & 0.874 & 0.746 & 0.828 & 0.816 \\
Tennessee (TN)   & 0.892 & 0.768 & 0.815 & 0.825 \\
\midrule
\textbf{Average} & \textbf{0.883} & \textbf{0.721} & \textbf{0.811} & \textbf{0.805} \\
\bottomrule
\end{tabular}
\end{table}

The GAN reliably reproduces marginal distributions while achieving strong manifold coverage and authenticity across states. The composite performance above 0.80 indicates high-quality synthetic scenario generation suitable for downstream resilience stress testing, consistent with benchmarks for precision–recall evaluation of generative models and tabular data synthesis \shortcite{kynkaanniemi2019improved,Xu2019CTGAN}.

\subsection{System Performance Evaluation}

\subsubsection{Baseline Transportation Performance}

We construct a representative batch of shipment orders covering the ten highest-frequency OD pairs in the TVA corridor, spanning a two-week planning horizon divided into daily time periods. The ten shipments originate from five cities and flow toward major inland hubs. Order volumes range from approximately 2 to 63 metric tons per shipment, with individual availability windows and delivery deadlines assigned within the two-week window. Mode-specific per-unit transportation costs and CO\textsubscript{2} emission factors are sourced from \cite{LIU2026105513}, and delay penalty rates are calibrated from \shortcite{gong2012assessing}.

Under nominal network conditions without disruptions, the intermodal routing model is solved to establish the simulation baseline. The baseline yields a transportation cost of \$4.85M, a penalty cost of \$2.16k, and a carbon emission cost of \$69.3k, for a total system cost of \$4.92M.

\subsubsection{Performance Under Historical Disruptions}

To assess system resilience against known events, we run simulations under three historically documented disruptions in the TVA region. For each event, affected nodes or waterway links are disabled during the disruption window and the optimization model is re-solved to obtain the system cost.

\rev{The disruption cases include Mississippi River Flooding in Spring 2023, reported by \citeN{noaa_mississippi_2023} as impairing barge operations along the Memphis--St.~Louis corridor; the Tennessee--Kentucky Tornado Outbreak in Spring 2024, documented by \citeN{wiki_tornado_2024} as damaging regional transportation infrastructure; and the Mississippi River Low Water Event in Late 2023, described by \citeN{princeton_navigation_2023} as restricting navigation drafts along the river system. The performance is listed in Table~\ref{tab:historical}.}

\begin{table}[!ht]
\centering
\caption{System performance under real historical disruptions.}
\begin{tabular}{lrrrr}
\toprule
Case & Transportation (\$) & Penalty (\$) & Carbon Cost (\$) & Total Cost (\$) \\
\midrule
Baseline 
& 4,848,566.69 
& 2,157.99 
& 69,303.77 
& 4,920,028.47 \\

Mississippi Flood 2023 
& 5,000,244.82 
& 19,667.30 
& 72,455.33 
& 5,092,367.45 \\

Tornado Outbreak 2024 
& 4,864,952.26 
& 13,294.37 
& 66,153.66 
& 4,944,400.29 \\

Low Water 2023 
& 4,848,566.69 
& 2,157.99 
& 69,303.77 
& 4,920,028.47 \\
\bottomrule
\end{tabular}
\label{tab:historical}
\end{table}

Relative to baseline conditions, the Mississippi Flood increases total cost by about 3.5\%, mainly due to higher penalty and rerouting costs, accompanied by increased carbon emissions from detours and modal shifts. The Tornado Outbreak leads to a modest cost increase of approximately 0.5\%, reflecting partial redundancy and alternative routes. In contrast, the Low Water scenario produces identical results to the baseline, suggesting that the affected river links do not constrain network feasibility under the modeled demand structure. Across these experiments, results reveal heterogeneous resilience: large-scale flooding generates substantial cost escalation, while localized disruptions are often absorbed by network flexibility.

\subsubsection{Stress Test Results}

To evaluate forward-looking resilience, we run simulation experiments under three GAN-generated disruption scenarios representing high-probability regional failure patterns inferred from the learned frequency distribution. These scenarios target critical river corridors and high-degree transportation nodes.

\begin{table}[!ht]
\centering
\caption{System performance under GAN-generated stress scenarios.}
\begin{tabular}{lrrrr}
\toprule
Scenario & Transportation (\$) & Penalty (\$) & Carbon Cost (\$) & Total Cost (\$) \\
\midrule
Tennessee River Flood
& 6,201,216.85 
& 4,065.34 
& 63,003.78 
& 6,268,285.96 \\

Ohio River Navigation Disruption 
& 4,848,566.69 
& 2,157.99 
& 69,303.77 
& 4,920,028.47 \\

Lower Mississippi Severe Storm 
& 4,864,952.26 
& 4,065.34 
& 66,153.66 
& 4,935,171.26 \\
\bottomrule
\end{tabular}
\label{tab:GAN}
\end{table}

Table~\ref{tab:GAN} shows that the Tennessee River flood scenario produces a substantial cost increase of approximately 27\% relative to baseline, indicating high structural vulnerability when multiple upstream river nodes are simultaneously disrupted. In contrast, the Ohio River navigation restriction produces no measurable cost impact, revealing strong redundancy along that corridor. The Lower Mississippi storm scenario results in a modest increase of approximately 0.3\%, suggesting localized but manageable rerouting adjustments.

We quantify systemic risk by weighting each scenario's cost outcome based on its probability, estimated from the GAN-generated frequency distribution. Let $C_i$ denote the total system cost under scenario $i$ and $\pi_i$ its associated probability weight; the probability-weighted expected system cost is \rev{\(\mathbb{E}[C] = \sum_i \pi_i C_i\)}.

\begin{table}[!ht]
\centering
\caption{Probability-weighted disruption scenarios.}
\label{tab:prob_scenarios}
\begin{tabular}{lccccccc}
\toprule
 & Baseline & West TN & Nashville & Paducah & East TN & North AL & Central TN \\
\midrule
Probability 
& 0.22 & 0.18 & 0.14 & 0.12 & 0.10 & 0.09 & 0.07 \\

Total Cost (\$) 
& 4,920,028 
& 4,933,168 
& 4,928,256 
& 4,985,378 
& 4,920,028 
& 5,936,198 
& 6,224,805 \\

Extra Loss (\$) 
& - 
& 13,140 
& 8,228 
& 65,350 
& 0 
& 1,016,169 
& 1,304,776 \\
\bottomrule
\end{tabular}
\end{table}

Table~\ref{tab:prob_scenarios} reports scenario probabilities, realized system costs, and the incremental loss above baseline for the seven highest-frequency disruption regions. The probability-weighted expected cost equals approximately \$5.11 million annually, representing an expected disruption premium of roughly \$190,000 above nominal operating conditions. Additionally, nearly the entire expected loss is concentrated in North Alabama and Central Tennessee, which individually generate cost increases exceeding 20\%.

This concentration effect reveals that systemic risk is not evenly distributed across the network. Most localized corridor or single-node failures produce negligible cost deviations and are absorbed by routing flexibility. In contrast, correlated multi-node inland disruptions create disproportionate economic impact, indicating structural vulnerability concentrated in specific regional clusters.

\subsubsection{Redundancy and Improvement Evaluation}

The probability-weighted results show that several disruption regions leave system cost unchanged, motivating a complementary single-node ablation experiment. We disable each major river port individually across separate simulation runs and record the realized total cost relative to the baseline, thereby isolating each node's structural contribution to network performance.

\begin{table}[!ht]
\centering
\caption{Single-node disruption impact.}
\label{tab:node_criticality}
\begin{tabular}{lrr}
\toprule
Node & Total Cost (\$) & Extra Loss (\$) \\
\midrule
Baseline (No Disruption) & 4,920,028 & -- \\

Port of Louisville & 4,920,028 & 0 \\
Port of Knoxville & 5,936,198 & 1,016,169 \\
Port of St.~Louis and East St.~Louis & 4,935,171 & 15,143 \\
Paducah--McCracken Riverport & 4,940,861 & 20,833 \\
Ports of Cincinnati--Northern Kentucky & 4,920,028 & 0 \\
The Port of Nashville & 4,928,256 & 8,228 \\
Port of Memphis & 4,942,546 & 22,518 \\
Port of Huntington Tri-State & 4,920,028 & 0 \\
\bottomrule
\end{tabular}
\end{table}

Table~\ref{tab:node_criticality} reveals pronounced heterogeneity in node importance. \rev{Disabling the {Port of Knoxville} increases total system cost by over \$1 million, indicating that this facility functions as a critical routing hub under the current freight demand configuration. Its removal forces substantial rerouting and modal substitution, producing cost escalations exceeding 20\%.}

In contrast, disabling the ports of Louisville, Cincinnati-Northern Kentucky, or Huntington, produces no measurable change in total cost. Under the modeled demand structure, these facilities are either redundant or not binding for optimal routing decisions. This does not imply that the nodes are operationally unimportant in general, but rather that they are not cost-critical for the specific OD flows considered.

Intermediate effects are observed for ports of Memphis, Paducah-McCracken, Nashville, and St.~Louis and East St.~Louis, which generate modest cost increases below \$25,000. These facilities contribute to routing flexibility but are partially substitutable within the broader network.

\section{CONCLUSION}
\label{sect:5}

This paper develops a data-driven stress-testing framework for evaluating resilience in intermodal freight networks, integrating GAN-based scenario generation. Applied to the TV corridor, historical disaster data are used to train a per-region GAN that produces correlated, multi-node disruption scenarios; each scenario is evaluated with the time-expanded intermodal routing model under disrupted network conditions. Repeating this process across the full scenario ensemble produces an empirical distribution of system cost outcomes, enabling both expected cost estimation and tail-risk analysis.

The results reveal strongly asymmetric vulnerability across disruption types. Historical events such as Mississippi River flooding increase total system cost by about 3.5\%, while most localized disruptions have negligible impact due to network redundancy. In contrast, GAN-generated compound inland disruptions increase costs by over 25\%, indicating that correlated multi-node failures drive disproportionate system-wide effects. Probability-weighted analysis estimates an expected system cost of approximately \$5.11 million, with a disruption premium of about \$190,000, largely concentrated in inland regions such as North Alabama and Central Tennessee. \rev{Node-level analysis further identifies the Port of Knoxville as a critical hub whose disruption results in a cost increase exceeding \$1 million.} These findings show that systemic risk is dominated by rare, correlated regional failures rather than frequent localized disruptions. The proposed framework provides a practical tool for identifying critical infrastructure, quantifying economic impacts, and prioritizing resilience investments under realistic disruption conditions.

This study has several limitations: The GAN relies on available historical data and may underrepresent extremely rare tail events. The simulation model captures cost impacts under capacity disruptions but does not incorporate stochastic congestion dynamics or cascading interdependencies across infrastructure systems. \rev{Future implementations can also test alternative regional partitions when they are more useful for planning.} Future work will extend the framework to incorporate dynamic traffic effects, richer compound hazard modeling, and adaptive decisions to enhance its applicability for forward-looking resilience planning.

\section*{ACKNOWLEDGMENTS}

The views and opinions expressed in this article are those of the authors and do not necessarily reflect the position of Amazon.com, Inc. or its affiliates. This work was conducted independently of the authors' roles at their respective institutions.

\footnotesize

\bibliographystyle{wsc}
\bibliography{references_revised}

\section*{AUTHOR BIOGRAPHIES}

\noindent {\bf \MakeUppercase{Xudong Wang}} is with the Department of Industrial and Systems Engineering (ISE) at the University of Tennessee, Knoxville. His research interests include generative AI and urban resilience systems. His email address is \email{xwang97@vols.utk.edu}.\\

\noindent {\bf \MakeUppercase{Mustafa Can Camur}} \rev{is with Amazon.} He received his Ph.D. in Industrial Engineering from Clemson University. His research interests include optimization and simulation. His email address is \email{mcamur@utk.edu}.\\

\noindent {\bf \rev{\MakeUppercase{Sabarna Choudhuri}}} \rev{is with Amazon.} He holds a Master’s degree in Electrical Engineering from Arizona State University. His Linkedin profile is available at \email{https://www.linkedin.com/in/sabarna-choudhuri/}. His e-mail address is \email{sabarna121@gmail.com}.\\

\noindent {\bf \MakeUppercase{Xueping Li}} is a Professor and Dan Doulet Faculty Fellow of ISE and the Director of iLab and co-Director of HITS Lab at UTK. He holds a Ph.D. from Arizona State University. He is an IISE Fellow. His e-mail address is \email{Xueping.Li@utk.edu
}. \\

\end{document}